\documentclass[11pt]{article}
\usepackage{a4wide,color}
\usepackage[centertags,leqno]{amsmath}
\usepackage{amssymb}
\usepackage[mathscr]{eucal}
\usepackage{epsf}

\usepackage{color}

\newtheorem{theorem}{Theorem}[section]
\newtheorem{prop}[theorem]{Proposition}
\newtheorem{lemma}[theorem]{Lemma}

\newtheorem{remark}[theorem]{Remark}
\newtheorem{cor}[theorem]{Corollary}

\newtheorem{definition}[theorem]{Definition}

\newcommand{\qed}{\rule{2mm}{2mm}}

\newcommand{\eqdef}{\stackrel{{\mathrm {def}}}{=}}

\newcommand{\proof}{{\em Proof. }}

\newcommand{\RR}{\mathbb{R}}

\numberwithin{equation}{section}

\begin{document}

\title{\bf
       Radial singular solutions for the $N-$ Laplace equation
       with exponential nonlinearities\\
\vspace{0.3cm}
}

\author{%
  Marius {\sc Ghergu}\\
          School of Mathematics and Statistics, \\
          University College Dublin,\\
          Belfield, Dublin 4, Ireland\\
{\tt marius.ghergu@ucd.ie}\\
and\\
Institute of Mathematics Simion Stoilow of the Romanian Academy,\\
21 Calea Grivitei Street,\\ 
010702 Bucharest, Romania\\
\vspace{0.1cm}
\and
        Jacques {\sc Giacomoni}\\
        LMAP (UMR E2S-UPPA CNRS 5142)\\
        Universit\'e de Pau et des Pays de l'Adour\\
        Avenue de l'Universit\'e,
        F-64013 Pau cedex, France\\
{\tt jacques.giacomoni@univ-pau.fr}\\
\vspace{0.1cm}
\and
        S. {\sc Prashanth} \\
        TIFR-Centre for Applicable Mathematics,\\
        Post Bag No. 6503, Sharadar Nagar,\\
        GKVK Post Office,\\
        Bangalore 560065, India\\
{\tt pras@math.tifrbng.res.in}\\
\vspace{0.1cm}
}

\maketitle

\vspace{-0.5cm}

\date{}
\noindent
\begin{abstract}
\baselineskip=14pt
In this paper, we consider radial distributional solutions of the quasilinear equation $-\Delta_N u=f(u)$ in the punctured open ball $ B_R\backslash\{0\}\subset \RR^N$, $N \geq 2$. We obtain  sharp conditions on the nonlinearity $f$ for extending such solutions to the whole domain $B_R$ by preserving the  regularity.  For a certain class of noninearity $f$ we obtain the existence of singular solutions and deduce upper and lower estimates on the growth rate near the singularity.
\end{abstract}

\vfill
\par\vspace*{0.5cm}
\noindent
\begin{tabular}{ll}
{\bf running head:}
& Radial Singular solutions for $N$-Laplace  equation\\
\end{tabular}

\par\vspace*{0.5cm}
\noindent
\begin{tabular}{ll}
{\bf Keywords:}
&  $N$-Laplace equation;\\
& removable singularity;\\
& singular solutions;\\
& behaviour near the singularity\\
\end{tabular}

\par\vspace*{0.5cm}
\noindent
\begin{tabular}{lll}
{\bf 2000 Mathematics Subject Classification:}
& Primary   & 35J62, 35B09;\\
& Secondary & 35B65 \\
\end{tabular}


\baselineskip=16pt

\section{Introduction}
\label{s:Intro}

Let $B_R \subset  \RR^N$, $N\geq 2$, be the open ball of radius $R$  centered at the origin in $\RR^N$. Consider the problem
\begin{eqnarray*}
( P_R^*)\qquad \left\{\begin{array}{ll}
 & - \Delta_N u
 =  f(u)
  \quad \mbox{ in }\,\ {\mathcal D}^\prime(B_R\setminus\{0\}) ,\\
&   u \geq 0 \quad \mbox{ in  }\, B_R\setminus\{0\},\\
 &  u \in L^\infty_{\rm loc}(\overline B_R\setminus\{0\})\,, \;\; u \text{ is radial}.
\end{array}\right.
 \end{eqnarray*}
Throughout the work, we assume  that $f$ satisfies:
\medskip

\noindent {\it 
$(f_1)$  $f\in C^{0,\theta}[0,\infty)$ for some exponent $\theta\in (0,1) $;} 
\smallskip

\noindent {\it $(f_2)$\;\; $\displaystyle \liminf_{t \to \infty} f(t)\; t^{1-N}>0; $}
\smallskip

\noindent {\it $(f_3)$  the map $t \mapsto f(t)+ \kappa t^{N-1}$ is non-decreasing for some $\kappa>0$.}
\smallskip

By a result of {\sc Brezis} and {\sc Lions} \cite{BrLi} (in the case $N=2$) and {\sc V\'eron} \cite[theorem 5.10, p.~283]{Ve1} (in the case $N>2$),  we have the following
\begin{theorem}\label{1}
 If $u $ solves $(P^*_R)$, then $u \in W^{1,p}(B_R)$ for any $1 \leq p<N$. Furthermore, 
 there exists some $\alpha \geq 0$ such that $u$ solves the following problem :
\begin{eqnarray*}
( P_{\alpha,R})\qquad \left\{\begin{array}
{ll}
 & - \Delta_N u
 =  f(u)+\alpha\delta_0
  \quad \mbox{ in }\,\ {\mathcal D}^\prime(B_R) ,\\
&f(u)\in L^\infty_{\rm loc}(\overline B_R\setminus\{0\})\cap L^1(B_R).
\end{array}\right.
 \end{eqnarray*}
\end{theorem}

\noindent The above result leads naturally to the following questions:
\medskip

\noindent {\it (Q1) Is there a sharp growth condition on $f$  that  determines whether a  solution of $(P^*_R)$ can or cannot be extended to a (distributional) solution to $(P_{0,R})$ ?}\medskip

\noindent {\it (Q2) If such an extension holds and the solution in the punctured domain $B_R\setminus\{0\}$ is smooth, is the extended solution equally smooth in $B_R$?}\medskip

\noindent {\it (Q3) If the extended solution blows up at the origin, what is its asymptotic blow-up rate?}
\medskip

When $N=2$ and $f$ has at most a polynomial growth, the first two questions $(Q1)$ and $(Q2)$ were discussed in detail in {\sc Brezis-Lions} \cite{BrLi} and {\sc Lions} \cite{Li-JDE}.
 For a corresponding discussion involving exponential growth nonlinearities in the case  $N=2$, we refer to {\sc Dhanya-Giacomoni-Prashanth} \cite{DhGiPr-AMS}.

In the quasilinear case $N \geq 3$, without being exhaustive,  we  mention the results of {\sc Guedda-V\'eron} \cite{GuVe}, {\sc Bidault-V\'eron} \cite{Bidault-Veron-arma} and {\sc Kichenassamy-V\'eron} \cite{KiVe}. For more  on the subject, we refer to the  survey {\sc V\'eron} \cite{Veron-handbook} and the book {\sc V\'eron} \cite{Ve1}. We cite the book {\sc Ghergu-Taliaferro} \cite{GhTa} for results concerning  isolated singularities for partial differential inequalities.

In the present paper, we extend the results contained in \cite{DhGiPr-AMS} to the radial quasilinear case  and obtain a complete answer to questions $(Q1)$ and $(Q2)$ above. We also provide a partial answer to $(Q3)$. 

\begin{definition}\label{def1.1}
We call $f$ a sub-exponential function if  there exists $\beta>0$ such that
\begin{equation}\label{eq1.1}
\displaystyle C:= \sup_{[0,
\infty)} f(t)e^{-\beta t} < \infty.
\end{equation}
 We call  $f$ to be super-exponential if it is not a
sub-exponential function.
\end{definition}
As a complete answer to question $(Q1)$, we show that if $f$ is sub-exponential then we can construct, by  the method of monotone iterations, a solution $U$ to $(P^*_R)$ which solves $(P_{\alpha,R})$ for some $\alpha>0$ (see Theorem~\ref{dirac-sol}). Conversely, we show in Lemma~\ref{13} that  any solution to $(P^*_R)$ extends to a solution to $(P_{0,R})$ if $f$ is super-exponential. In contrast to this result we show that the Dirac mass in $(P_{\alpha,R})$ is not in general removable for sub-exponential functions $f$.

Regarding question $(Q2)$ above, for a super-exponential $f$, we construct examples of radial solutions to $(P_{0,R})$ which blows up only at the origin (see Lemma~\ref{mu-lemma} and Theorem~\ref{critical-behaviour-f}). That is, such solutions are smooth in $B_R\setminus\{0\}$, but not in $B_R$. We also show that if $f$ is sub-exponential, then any solution $u$ of $(P^*_R)$ that extends to a solution of $(P_{0,R})$ is regular, say in $C_{loc}^{1,\theta}$ (see Theorem~\ref{Brezis-Merle-thm}). 

The question $(Q3)$ is partially answered in Lemma~\ref{19}. Utilising the asymptotic analysis of {\sc Atkinson} and {\sc Peletier} (see \cite{AtPe}), for the super-exponential nonlinearities $f$  we derive  an upper bound for singular solutions (see Lemma \ref{19}) and consequently obtain the following limiting behaviour for any  solution $u$ of $(P^*_R)$ :

\begin{equation*}
|x|^N f(u(x) ) \to 0   \; \; \text{as }\; \; |x| \to 0.
\end{equation*}

We remark that this result is new even for the semilinear case.  Although we do not obtain a pointwise  lower bound, we derive in  Lemma \ref{vic} an integral bound for the behavior of the solution around the isolated singularity at the origin; see also Corollaries \ref{coo} and \ref{bb}.
The accurate asymptotic behaviour of singular solutions is still an open question. Nevertheless, Corollaries \ref{hon} and \ref{meh} give alternatives similar to the ones available in higher dimensions  for supercritical nonlinearities (see for instance Theorem 5.13 in \cite{Ve1}).  

The  semilinear case $N=2$ with exponential-type nonlinearities in dimensions higher  than 2 is dealt with in the recent work of {\sc Kikuchi} and {\sc Wei} \cite{KiWe}. By using  the Emden-Fowler transformation and a clever analysis of the perturbation term in the  asymptotic profile, the authors prove the existence and the precise asymptotic behaviour of singular solutions. In  the two dimensional case, this approach  as well as the one using Harnack type inequalities as in \cite{Bidault-Veron-arma}, fails since the coefficients of the resulting equation do not posses the integrability in the right spaces.

\section{ Some preliminary results}\label{1.5}
In this section we collect some useful results to our approach.
The following first observation is in order:
 \begin{prop}\label{p2}
 If $u$ solves $(P^*_R)$, then $u$ is radially non-increasing in $B_R\setminus\{0\}$.
 \end{prop}
 \proof In the radial variable, $u$ solves the ODE:
 $$
 -r^{1-N}\Big( r^{N-1} |u^\prime|^{N-2} u^{\prime} \Big ) ^{\prime} = f(u) \;\;\mbox{ in } \;\; (0,R).
 $$
 Noting that $f(u) \geq 0$, it readily follows that the map $r \mapsto r u^{\prime}(r)$ is non-increasing in $(0, R)$. Let
 \begin{equation}\label{0.5}
 \ell:= \lim_{r \to 0^+} r u^{\prime}(r).
 \end{equation}
 If $\ell>0$, by integrating \eqref{0.5} we obtain that $u(r) \to - \infty$ as $r \to 0^+$, a contradiction to the nonnegativity of $u$ near $0$.  Hence $\ell \leq 0$ and the result follows.
 \hfill \qed 

 First, we note the following simple consequence of definition \ref{def1.1}. 
\begin{prop} \label{p1} If $f$ is a super-exponential nonlinearity, then
$$\int_0^{\infty} f(t) e^{-\epsilon t} dt = \infty \quad \mbox{ for any } \epsilon>0.$$
\end{prop}
\proof From \eqref{eq1.1} we have that for any $\epsilon>0$  and any positive integer $n$ 
$$
\sup_{[n,\infty)} f(t)e^{-\epsilon t} = \infty.
$$
Therefore, we may find $t_n \geq n$ such that 
$$
f(t_n) \geq n e^{\epsilon t_n} . \quad  
$$
Thus,
$$ \int_{t_n}^{\infty}( f(t) + \kappa t^{N-1})e^{-\epsilon t} dt  \geq  \int_{t_n}^{\infty} (f(t_n) + \kappa t_n^{N-1})e^{-\epsilon t} dt \geq \frac{n }{\epsilon}.
$$
The assertion now follows. \hfill \qed

Our next result establishes the connection between distributional and entropy solutions to $(P_R^*)$.
We recall that for $g \in L^1(\Omega)$, a function $u \in L^1_{loc}(\Omega)$ is an entropy solution of 
\begin{equation} \label{e0.01}
-\Delta_N u = g \;\; \mbox{ in } \;\; \Omega,  \quad 
u = 0 \;\; \mbox{ on } \;\;\partial \Omega,
 \end{equation}
 if $T_k(u) \in W^{1,N}_{0}(\Omega)$ for any $k \in \mathbb{R}^+$,   $|\nabla u| \in L^{N-1}_{loc}(\Omega)$ and
\begin{equation}\label{0.1}
\int_\Omega |\nabla u|^{N-2} \nabla u \cdot \nabla T_k(u-\phi) \leq \int_\Omega  T_k(u-\phi) g \quad \mbox{ for all } \phi \in C_0^{\infty}(\Omega),  k \in \mathbb{R}^+.
\end{equation}
Here, for $k \in \mathbb{R}^+,$ $T_k$ denotes the truncation map given by
$T_k(t) =t$ for  $t \in [-k,k]$ and $T_k(t)=k \,{\rm sgn}(t)$  for $|t| \geq k$.

\begin{prop}\label{p0.5}
Let  $u \in C^1_{loc}(\overline B_R\setminus\{0\}) \cap L^1_{loc}(B_R), f \in L^1(B_R)$ be  both radial functions and  \\ $|\nabla u|^{N-1} \in L^{p}(B_R)$ for some $p>1$. If $u$ is a distributional solution to $-\Delta_N u = f$ in $B_R$, then $u-u(R)$ is an entropy solution to the same problem in $B_R$ with homogeneous Dirichlet data on $\partial B_R$.
\end{prop}
\proof  Since $u$ is a distributional solution we have
$u\in L^1_{\rm loc}(B_R)$ has "zero trace" on $\partial\Omega$, $\vert \nabla u\vert\in L^{N-1}_{\rm loc}(B_R)$ and  
\begin{equation}\label{0.2}
\int_{B_R}\vert\nabla u\vert^{N-2}\nabla u\cdot\nabla\phi =\int_{B_R} \phi \, g  \quad \mbox{ for all } \phi\in C^\infty_0(B_R).
\end{equation}

We may fix $p$ close to 1 so that  $W^{1,p^{\prime}}_0(B_R) \hookrightarrow C_0(\overline{B}_R)$. We first claim that if $|\nabla u|^{N-1} \in L^{p}(B_R)$ for some $p>1$ and $ f \in L^1(B_R)$,  we can enlarge the class of test functions $\phi$  in  \eqref{0.2}  to $W^{1,p^{\prime}}_0(B_R)$. To see this, we take a sequence $\{\phi_n\} \subset C_0^{\infty}(B_R)$ converging to $\phi $ in the $W^{1,p^{\prime}}_0(B_R)$ norm. Then $\phi_n \to \phi$ in $C_0(\overline{B}_R)$ and the claim follows by writing \eqref{0.2} for each $\phi_n$ and passing to the limit $n \to \infty$.

 Let $v:= u-u(R)$.  
  Given $\phi \in C_0^{\infty}(B_R)$, we  have that $T_k(v-\phi) \in W^{1,\infty}_0(B_R)$ for every $k \in \mathbb{R}^+$. Hence by the above observations,
  \begin{eqnarray}
  \int_{B_R} |\nabla v|^{N-2} \nabla v \cdot \nabla T_k(v-\phi) &= &   \int_{B_R} |\nabla u|^{N-2} \nabla u  \cdot \nabla T_k(v-\phi)  \notag \\
 &=& \int_{B_R} f T_k (v-\phi)  \quad \forall \phi \in C_0^{\infty}(\Omega), \; \forall k \in \mathbb{R}^+. \notag
 \end{eqnarray}
 Thus, $v$ satisfies \eqref{0.1}. \hfill \qed

Finally, we recall a version of {\sc Brezis-Merle} \cite{BrMe} result for  the $N$-Laplace problem. 

\begin{prop}\label{p0.6} (see \cite[Theorem 1.6]{AgPe})
Let $u$ be the entropy solution of the problem \eqref{e0.01} where $g \in L^1(\Omega)$. Then the following inequality holds for any $\delta \in (0, N \omega_N^{\frac{1}{N-1}})$:
\begin{equation}\label{e0.55}
\int_\Omega exp \Bigg(\frac{\big(N \omega_N^{\frac{1}{N-1}} - \delta\big)|u|}{\| g\|_{1,\Omega}^{\frac{1}{N-1}}} \Bigg) \leq \frac{N\omega_N^{\frac{1}{N-1}}}{\delta}|\Omega|.
\end{equation}
Here, $\| \cdot \|_{1, \Omega}$ denotes the $L^1(\Omega)$ norm and $|\Omega|$ the N-dimensional Lebesgue measure of $\Omega$ and $\omega_N$ is the volume of the unit sphere in $\RR^N$.
\end{prop}
\begin{cor}\label{c1} (see \cite[Corollary 1.7]{AgPe})
Let $u$, $g$ be as above. Then, $e^{|u |}\in L^p(\Omega)$ for any $p \geq 1$.
\end{cor}
\section{ Dirac mass solution for sub-exponential $f$}\label{section2}
 In this section, we show the following:
\begin{theorem}\label{dirac-sol}
Let $f$ be a sub-exponential function  and  $\beta, C$ be  given by \eqref{eq1.1}.
 Then there exist $\alpha, R_*>0$ depending on $\beta, C$ such that 
  $(P_{\alpha,R_*})$  admits a distribution solution. 
   \end{theorem}
Before proving the theorem, we construct appropriate sub and super solutions in the lemmas below.
\begin{lemma}\label{upper-sol}
Let $f$ be a sub-exponential function with $\beta, C$ be given by  \eqref{eq1.1}. Define 
\begin{equation}\label{4}
v_{\beta,C}=\frac{1-N}{\beta}
\Big [ \log r + \frac{N}{N-1} \log (1+ C \beta^{N-1} r) \Big], \;\; r>0.
\end{equation}
Then, the following pointwise inequality holds :
\begin{equation}\label{5}
-\Delta_N v_{\beta,C} \geq f(v_{\beta,C}) \;\; \text{ in } \RR^N \setminus \{0\}.
\end{equation}
\end{lemma}
\proof Let 
$$
\mu := \frac{N}{N-1}, \; M:= C \beta^{N-1}.
$$
By a straightforward calculation,
\begin{eqnarray}
-\Delta_N  v_{\beta,C} (r) &=&  r^{1-N}(-r^{N-1} |v_{\beta}^{\prime}(r)|^{N-2} v_{\beta}^{\prime})^\prime \notag \\
 &=& \Big(\frac{\mu M (N-1)^N}{\beta^{N-1}}\Big) \Big [\frac{1}{r^{N-1}(1+M r )^2} \Big] \Big[ \frac{1+(\mu+M r)}{1+M r} \Big]^{N-2}. \label{6}
\end{eqnarray}
We recall the growth condition on $f$ in definition \ref{def1.1} and use \eqref{4} to find that
\begin{equation}\label{7}
f(v_{\beta,C}) \leq C e^{\beta v_{\beta,C}} = \frac{C}{r^{N-1}(1+M r)^N}.
\end{equation}
Therefore, from \eqref{6} and \eqref{7} we obtain that
\begin{equation} -\Delta_N v_{\beta,C} \geq f(v_{\beta,C}) \;\; \text{ in } \RR^N \setminus \{0\}. \end{equation}   \hfill \qed
\begin{remark}\label{r1}
Note that $v_{\beta,C}$ is a strictly decreasing function on $(0,\infty)$ with $$lim_{r \to 0^+} \; v_{\beta,C} = \infty \;\;\text{ and }\;\; lim_{r \to \infty} \; v_{\beta,C} = -\infty.$$
\end{remark}
\begin{definition}\label{d1}
Given $\beta, C>0$ define $R_*:=R_*(\beta,C)$ to be the unique zero of $v_{\beta,C}$.
\end{definition}
\begin{lemma}\label{lower-sol}
Let $f, \beta, C$ be as in the previous lemma. For $R>0$ define
\begin{eqnarray}\label{sub-sol}
w_{\beta,C, R}(r) &:= & \frac{1-N}{\beta} \Big [ \log r + \frac{N}{N-1} \log (1+ C \beta^{N-1} R) \Big], \;\; r>0.\label{8}
\end{eqnarray}
Then
\begin{eqnarray}
w_{\beta,C,R} &\leq & v_{\beta,C} \;\text{ in }\; B_R, \label{10} \\
w_{\beta,C,R_*} &\geq& 0 \;\text{ in }\; B_{R_*}, \label{10.1} \\
w_{\beta, C,R}  (R)&=& v_{\beta,C} (R), \label{9} \\
w_{\beta,C,R}, \; v_{\beta,C} &\in& C^\infty_{loc}(\RR^N \setminus \{0\}), \label{11}\\
-\Delta_N w_{\beta,C,R} &=&\omega_{N}^{\frac{1}{N-1}}\Big( \frac{N-1}{\beta}\Big) \delta_0  \;\text{ in }\; {\mathcal D}^\prime(\RR^N) \label{12}.
\end{eqnarray}
\end{lemma}
The proof of Lemma \ref{sub-sol} follows by straightforward calculations. 

\noindent Now we are ready to proceed with the proof of  Theorem~\ref{dirac-sol}. Let
$\beta, C>0$ be given by  \eqref{eq1.1} and let  $v_{\beta,C}, w_{\beta,C,R}$ be the corresponding functions as in \eqref{4} and \eqref{8}. 
Choose $R_*$ as in definition \ref{d1}. Given $0<\epsilon<R_*$, define the annular region
$$
A_{\epsilon,R_*}:= B_{R_*} \setminus B_\epsilon.
$$
We now set up the following iteration procedure (solved in the $W^{1,N}$-weak sense):
\begin{eqnarray*}
(I_n)\qquad \left\{\begin{array}
{ll}
& u_0 := w_{\beta,C,R_*}; \\
& - \Delta_N u_n +  \kappa |u_n|^{N-2} u_n
 =  f(u_{n-1}) + \kappa |u_{n-1}|^{N-2} u_{n-1}
  \; \mbox{ in }\, A_{\epsilon,R_*},\\
&u_n = w_{\beta,C,R_*} \; \text{ on }\; \partial A_{\epsilon,R_*}; \\
&  u_n \in  C^{1,\theta} (\overline{A}_{\epsilon,R_*}) \; \text{for any }\; \theta \in (0,1) \; \text{ and } u_n \; \text{ radial}.
\end{array}\right.
 \end{eqnarray*}
It is standard to see that the quasilinear equation in $(I_n)$ is solvable in the $W^{1,N}-$weak sense and by H\"older regularity results (see for instance  {\sc Tolksdorf} \cite{To} and {\sc Lieberman} \cite{Li}) the required regularity of the solution is obtained.


Using the hypothesis $(f_3)$, \eqref{5}, \eqref{10.1} and an induction argument, we obtain that
\begin{equation} 
 u_0 := w_{\beta,C,R_*} \leq u_n \leq u_{n+1} \leq v_{\beta,C}  \;\text{ in } \; A_{\epsilon,R_*} \; \text{ for all } \; n. \notag
\end{equation}
 Since $v_{\beta,C}$ is bounded in $ A_{\epsilon,R_*} $ for a fixed $\epsilon>0$, from the above pointwise estimates  we can pass to the uniform $C^{1,\theta}(\overline{A}_{\epsilon,R_*})$ estimates for the sequence $\{u_n\}$ by using the classical quasilinear regularity results (see for instance {\sc Tolksdorf} \cite{To} and {\sc Lieberman} \cite{Li}). Therefore, we obtain a function $u_{\epsilon,R_*}$ such that 
$$
u_n \to u_{\epsilon,R_*} \; \text{ in } C^{1}(\overline{A}_{\epsilon,R_*})
$$
and  $ u_{\epsilon,R_*} $ solves (in the $W^{1,N}$-weak sense) the following problem:
\begin{eqnarray*}
\left\{\begin{array}
{ll}
& - \Delta_N  u_{\epsilon, R_*} 
 =  f( u_{\epsilon, R_*} ) ,  w_{\beta, C , R_*} \leq  u_{\epsilon, R_*}  \leq v_{\beta,C} 
  \; \mbox{ in }\, A_{\epsilon, R_*},\\
&u_{\epsilon, R_*}  = \; w_{\beta, C , R_*}  \text{ on }\; \partial A_{\epsilon, R_*},\\
& u_{\epsilon, R_*}  \in  C^{1,\theta} (\overline{A}_{\epsilon, R_*})   \; \text{ and is radial}.
\end{array}\right.
 \end{eqnarray*}
 Choosing a positive sequence $\{\epsilon_n\}$ tending to $0$ and noting that the corresponding sequence of solutions $\{u_{\epsilon_n,R_*} \}$ is relatively compact in $C^1_{loc}(\overline{B}_{R_*} \setminus \{0\})$, we obtain a distributional solution $U$ of
 \begin{eqnarray*}
(P^*_{R_*}) \qquad \left\{\begin{array}
{ll}
& - \Delta_N  U 
 =  f( U ) , \quad  w_{\beta, C , R_*} \leq  U  \leq  v_{\beta,C} 
  \; \mbox{ in }\, B_{R_*}^*,\\
&U  =   w_{\beta, C , R_*}  \; \text{ on }\; |x|= R_*, \\
& U  \in  C^{1,\theta} (\overline B_{R^*}\setminus\{0\})   \; \text{ and is radial}.
\end{array}\right.
 \end{eqnarray*}
 Note that 
 $$
 U (r) \sim \Big(\frac{1-N}{\beta}\Big) \log r \;\text{ as }\; r \to 0^+ , 
\mbox{ and }  U, f(U)\in L^\infty_{\rm loc}(\RR^N \setminus \{0\}).
$$
 By Theorem \ref{1}  we obtain that 
 $$
 f(U) \in L^1_{loc}(B_{R_*}), \; |\nabla U | \in L^p(B_{R_*}) \;\; \mbox{ for all } \; 1\leq p < N,
 $$
 and for some $\alpha \geq 0$, 
  $$
  - \Delta_N  U =  f(U)+  \alpha \delta_0 \;\text{ in }\; {\mathcal D}^\prime(B_{R_*}).
 $$
 If $\alpha=0$, from Propositions \ref{p0.5}-\ref{p0.6}, Corollary \ref{c1} and  Corollary 2.2 in \cite{AgPe} we obtain  that $U$ is bounded near the origin, contradiction. Hence, necessarily $\alpha > 0$.
\hfill \qed
%

%
\section{Removable Singularity}\label{section3}

In this section we show that if $f$ is super-exponential then distributional solutions of $(P^*_R)$ can be extended to (distributional) solutions of $(P_{0,R})$.

\begin{lemma}\label{13}
Let $f$ be super-exponential  and let $u$  be a solution of  $(P^*_R)$, and hence of $(P_{\alpha,R})$ for some $\alpha \geq 0$ (see Theorem \ref{1}).
 Then, necessarily $\alpha=0$.
\end{lemma}
\proof Let $0<\eta<R$ be small. Choose a  nonnegative radial test function $\phi \in C_0^{\infty}(B_\eta)$ with $\phi(0)=1$ and $\max_{[0,\eta]}|\phi^\prime| \leq 1/\eta.$ We then obtain
\begin{equation}\label{14}
  \int_0^\eta r^{N-1}|u^\prime|^{N-1} \phi^\prime(r) dr = \int_0^\eta r^{N-1}f(u(r)) \phi(r) dr + \alpha \phi(0).
\end{equation}
We estimate
\begin{eqnarray}
\text{ L.H.S. of }\; \eqref{14} &=& O\Big(\max_{[0,\eta]}r^{N-1}|u^\prime|^{N-1}  \int_0^\eta |\phi^\prime| \Big) \notag \\
&=& O \Big (\max_{[0,\eta]}r^{N-1}|u^\prime|^{N-1}\Big) \label{15}. 
\end{eqnarray}
Since 
$$
\int_0^\eta r^{N-1}f(u(r)) \phi(r) dr \to 0 \;\text{ as }\;\eta \to 0^+,
$$
from \eqref{14} and \eqref{15} we will obtain $\alpha=0$ if we can show
\begin{equation}\label{16}
\max_{[0,\eta]}r^{N-1}|u^\prime|^{N-1} \to 0 \;\text{ as }\; \eta \to 0^+.
\end{equation}
To this aim, we utilise the Edem-Fowler transformation given by
\begin{equation} \label{17}
t:=N\log(N/r), \;\; 0< r < R ; \; \; y(t) := u(r), \;\;  N\log(N/R)<t<\infty.
\end{equation}
Consequently, 
$$
y^\prime (t) = -\frac{r}{N}u^\prime (r).
$$
It can then be easily checked that if $u$ is a radial solution of  $-\Delta_N u = f(u)$ in $B_R\setminus\{0\}$, then $y$ solves the following Emden-Fowler type ODE:
\begin{equation} \label{18}
-\Big( |y^\prime|^{N-2}y^\prime \Big)^\prime = e^{-t} f(y) \;\text{ in }\; (N\log(N/R), \infty).
\end{equation}

Note that since $u$ is non-increasing (by Proposition \ref{p2}), one has $y^\prime \geq 0$.
Therefore, \eqref{16} holds iff $y^\prime(t) \to 0$ as $ t \to \infty$. We note from \eqref{18} that
$$
-\Big( |y^\prime|^{N-2}y^\prime \Big)^\prime = e^{-t} f(y) \geq 0 \;\text{ in }\; (N\log(N/R), \infty).
$$
This immediately implies that $y^\prime$ is a decreasing function on $(N\log(N/R), \infty).$ Let us denote
$$
\ell:= \lim_{t \to \infty} y^\prime(t).
$$

Noting that $y^\prime \geq \ell$ in $(N\log(N/R), \infty)$, we obtain  
$$
y(t) \geq y(N\log(N/R)) + \ell(t-N\log(N/R)) \; \quad \mbox{ for all } t \geq N\log(N/R).
$$
Hence, from the above inequality and the assumptions $(f_2) - (f_3)$ we find
\begin{eqnarray*}
\infty > \int_{B_{R/2}} \big( f(u) + \kappa u^{N-1}\big)dx  &=& \int_{N\log(2N/R)}^{\infty} \Big(f(y(t))+ \kappa y(t)^{N-1} \Big) e^{-t} dt  \notag \\
&\geq&  \int_{N\log(2N/R)}^{\infty} f(y(N\log(N/R)) + \ell (t-N\log(N/R))))e^{-t} dt,  \notag
\end{eqnarray*}
which contradicts Proposition \ref{p1} if $\ell>0$. Therefore, $\ell=0$ and hence $\alpha=0$. 
\hfill \qed
%
\section{Existence of singular solutions}\label{section4}
In this section we answer question $(Q2)$. We first show that when $f$ is sub-exponential, then any distributional solution of $(P_{0,R})$ is regular.
\begin{theorem}\label{Brezis-Merle-thm}
Let $f$ be a sub-exponential nonlinearity. Then any solution to $(P_{0,R})$ is regular in $B_R$.
\end{theorem}
\proof 
Let $u$ be a solution to $(P_{0,R})$. Then from propositions \ref{p0.5}- \ref{p0.6} and corollary \ref{c1}, we have that $e^{u}$ (and hence $f(u)$) belongs to $L^p(B_1)$ for all $p \geq 1$.  Thus, applying corollary 2.2 in \cite{AgPe} (see also related results in {\sc Ioku} \cite{Io} and {\sc Boccardo-Peral-Vazquez} \cite{BoPeVa}) we obtain that $u\in W^{1,N}_0(B_1)\cap L^\infty(\Omega)$. Hence from H\"older regularity results in \cite{Li}, we obtain that $u\in C^{1,\theta}(\overline{\Omega})\cap C_0(\overline{\Omega})$ for all $0< \theta<1$.
\hfill \qed
%
%

In the following two results, we construct solutions to $(P_{0,R})$ which blow up only at the origin for some classes of  super-exponential nonlinearities. 
\begin{lemma}\label{mu-lemma}
Given any $\mu>1$, there exists a  super-exponential nonlinearity $f$ satisfying :
$$\displaystyle\lim_{t\to \infty}f(t)e^{-t^\mu}=0, \;\; \displaystyle\lim_{t\to \infty}f(t)e^{-t^{\mu-\epsilon}}=\infty \;\; \forall\epsilon>0$$
 such that the corresponding problem $(P_{0,\frac{1}{2}})$ admits a radial solution that blows up only at the origin.
\end{lemma}
\proof Given $\mu>1$, define 

$$
f(t) := \quad  \Big(\frac{N}{\mu}\Big)^{N}(\mu-1)(N-1)  \left\{\begin{array}
{ll}
\; (N \log 2)^{\frac{(1-\mu)(N-1)-\mu}{\mu}} \; e^{(N \log 2)},  \;\; & 0 \leq t \leq (N \log 2)^{\frac{1}{\mu}}\\
 \; t^{(1-\mu)(N-1)-\mu}\; e^{t^\mu}, \;\; & t> (N \log 2)^{\frac{1}{\mu}}.
\end{array}\right.
$$

It can be easily checked that 

$$u(x) \eqdef \left(\log\frac{1}{\vert x\vert^N}\right)^{\frac{1}{\mu}}, \;\; x \neq 0$$ 

solves $(P^*_{\frac{1}{2}})$. Appealing to Theorem \ref{1} and Lemma \ref{13} we obtain that
$u$  solves $(P_{0,\frac{1}{2}})$. It can be also directly checked that $f(u) \in L^1(B_{\frac{1}{2}})$. %
\hfill \qed

In the next result, we exhibit a class of super-exponential  nonlinearities whose growth rate at infinity is of critical type (in the sense of Trudinger-Moser).
Although we state the result for nonlinearities in a general form, one can check that the model class of nonlinearities 
$$f_0(t):= t^{-\alpha} (1+t)^m e^{t^{\frac{N}{N-1}}-t^\beta}, \quad \alpha>0, \; m \geq 0, \; \frac{1}{N-1} < \beta < \frac{N}{N-1} $$
satisfy the required assumptions in (i).
\begin{theorem}\label{critical-behaviour-f}
(i) Let $f_0:(0,\infty)\to (0,\infty)$ be a smooth (at least $C^3$) nonlinearity  which has the form $f_0(t)=h(t)e^{t^{\frac{N}{N-1}}}$ where $h(t)$ satisifies assumptions (A1)-(A5) in  \cite{GiPrSr-asy}. Then there exists $R>0$ such that $-\Delta_N u = f_0(u)$  admits a radial positive distributional solution $z^*$ in $B_R$ blowing up only at the origin and such that  $z^* \not \in W^{1,N}_{loc}(B_R)$.\\
(ii) Let $f(t) :=  f_0(\frac{R}{2}) \chi_{[0, \frac{R}{2}]}(t) + f_0(t) \chi_{[\frac{R}{2}, \infty)}(t)$. Then  $(P_{0,\frac{R}{2}})$ admits a radial solution that blows up only at the origin.
\end{theorem}
\proof
(i) We  use  the transformation \eqref{17} in order to cast the differential equation $-\Delta_N u= f(u)$ in $B_R$ (for a positive solution)  into the following equivalent form :
%
%
%
\begin{equation}\label{ind}
\begin{array}{cllll}\left.
\begin{array}{rllll}
-(\vert y'\vert^{N-2}y')'& = & e^{-t}f_0(y), \\
y  \geq  0 \;, \quad y  &\in& L^{\infty}_{loc}\\
\end{array}\right\} \;\; \text{in} \;(T,\infty), \; T \eqdef 2 \log(2/R); \\
\hspace{-2in} \int_T^\infty f_0(y(t))e^{-t} dt<  \infty. 
\end{array}
\end{equation}
For our purposes, it is more convenient to consider the following ``shooting from infinity"  problem depending upon a parameter $\gamma>0:$
$$ (S_\gamma) \qquad \left \{\begin{array}{ccccc}
-(\vert y'\vert^{N-2}y')^\prime=e^{-t}f_0(y),\\
y(\infty)=\gamma,\; y'(\infty)=0.\\
\end{array} \right .
$$
Let $y(\cdot,\gamma)$ denote the unique solution of $(S_\gamma)$.  We see that $y(\cdot,\gamma)$ is a strictly concave function as long as it is nonnegative. Therefore, there exists a first zero of the solution $y(\cdot,\gamma)$ denoted by  $T(\gamma)$.  In addition, the map $\gamma\to T(\gamma)$ is continuous (see \cite[Lemma 3.1]{GiPrSr-asy}). We divide the proof into three steps:\\

{\it Step 1 : Given a sequence $\gamma_n\to\infty$ as $n\to\infty$, there exist a sequence $\{r_n\}$ of positive numbers with $\displaystyle\liminf_{n\to\infty}r_n>0$ and a sequence of  nonnegative radial solutions $\{u_n\} \subset L^{\infty}(B_{r_n})$ of $-\Delta_N u_n = f_0(u_n)$  in ${\cal D}^\prime(B_{R_n})$ with $u_n(0)=\gamma_n$.}\\

Given such $\{\gamma_n\}$, let $y_n=y(\cdot,\gamma_n)$ denote the solution of $(S_{\gamma_n})$. By \cite[proposition 4.2, p.~12]{GiPrSr-asy}, we obtain that $T^*=\displaystyle\limsup_{n\to\infty}T(\gamma_n)<\infty$. Up to a subsequence of $\{\gamma_n\}$, we can assume that $T^*=\displaystyle\lim_{n\to\infty}T(\gamma_n)$. We  fix this subsequence of $\{y_n\}$, which we still call $\{y_n\}$. Furthermore,  from the asymptotic behaviour of $f_0$ at $\infty$ it is not difficult  to show that $T^*>-\infty$ (see also \cite[lemma~4.1]{GiPrSr-JDE}).

Let us use the change of variable in \eqref{17} as follows:
$$r_n=Ne^{-T(\gamma_n)/N}, \;\; u_n(x)=y_n(N\log(N/\vert x\vert), \;\; x\in B_{r_n}\setminus\{0\}.$$ 
Then, we see that $u_n$ solves the required equation in $B_{r_n}$  with $u_n(0)=\gamma_n$ and $\displaystyle\liminf_{n\to\infty} r_n>0$.\\

{\it Step 2 : Let $\{\gamma_n\}$, $\{y_n\}$ be as in Step 1. Extend $y_n$ to $[T^*-1, T(\gamma_n))$ by $0$. Then the extended sequence (still denoted as $\{y_n\}$) is uniformly bounded on compact subsets of $[T^*-1,\infty)$.}\\

Denote $g(t)=\log(f_0(t))$. We define the following energy functional associated to $(S_{\gamma_n})$:
\begin{equation*}
E_n(t)=(y_n')^{N-1}(t)-\frac{N-1}{N}(y_n'(t))^Ng'(y_n(t))-e^{g(y_n(t))-t},\quad t>T(\gamma_n).
\end{equation*}
 Given $s_0>0$, let $t_0(\gamma_n)>T(\gamma_n)$ be  defined by the relation $y_n(t_0(\gamma_n))=s_0>0$. We may choose $s_0$ large enough such that $g',\, g''>0$ in $[s_0,\infty)$. Therefore, since $y_n$ is strictly increasing and $g$ is convex in $[t_0(\gamma_n),\infty)$,
\begin{equation*}
E_n'(t)=-\frac{N-1}{N}(y_n'(t))^{N+1}g''(y_n(t)) \leq 0, \quad\mbox{ for all } t\geq t_0(\gamma_n).
\end{equation*}
 Since $\displaystyle\lim_{t\to\infty}E_n(t)=0$, we obtain that $E_n$ is a nonnegative function on $[t_0(\gamma_n),\infty)$. This immediately implies that
\begin{equation}\label{eq3.1}
y_n'(t)g'(y_n (t))\leq \frac{N}{N-1},\quad \forall t\geq t_0(\gamma_n).
\end{equation}
Now, integrating the ODE in $(S_{\gamma_n})$, we have 
\begin{equation*}
(y_n^\prime(t_0))^{N-1}=\int_{t_0(\gamma_n)}^\infty f_0(y_n(t))e^{-t}\,\mathrm{d}t.
\end{equation*}
Therefore, from \eqref{eq3.1} and recalling that $y_n(t_0(\gamma_n))= s_0$, we obtain
\begin{equation}\label{eq3.2}
\displaystyle\sup_{n}\int_{t_0(\gamma_n)}^\infty f_0(y_n(t))e^{-t}\,\mathrm{d}t<\infty.
\end{equation}
Let $[a,b]\subset [T^*-1,\infty)$. Define $A=\{n\,:\, t_0(\gamma_n)>b\}$, $B=\{n\,:\, t_0(\gamma_n)\leq b\}$. We note that $\displaystyle\sup_n y_n(b)<\infty$. Otherwise, there must be a subsequence of $\{y_n(b)\}$ that tends to $\infty$ and hence by monotonicity of $y_n$, this subsequence converges uniformly to $\infty$ in $[b,b+1]$ violating \eqref{eq3.2}.

Again by the monotonicity of $y_n$, we have $\displaystyle\sup_{[a,b]}y_n\leq y_n(b)$. Therefore,
\begin{equation*}
\displaystyle\sup_{n\in A}\displaystyle\sup_{t\in [a,b]}y_n(t)\leq \displaystyle\sup_{n\in A}y_n(b)\leq \displaystyle\sup_n y_n(t_0(\gamma_n))=s_0.
\end{equation*}
Similarly, 
$$\displaystyle\sup_{n\in B}\displaystyle\sup_{[a,b]}y_n(t)\leq \displaystyle\sup_{n\in B}y_n(b)<\infty,
$$ 
which completes the proof in Step 2.\\

{\it Step 3:  Constructing the singular solution.} \\

From Step 2 and the fact that $y_n$ solves the ODE in \eqref{ind} we obtain a subsequence of $\{y_n\}$, which we denote again by $\{y_n\}$, such that   $\{y_n\}$ and $\{y_n'\}$ are uniformly convergent  in any compact sub-interval of $(T^*,\infty)$. By using a diagonalisation process, we can obtain a subsequence of $\{y_n\}$, which we will  denote by $\{y_n\}$ again, and a positive, continuous nondecreasing function $y^*$ on $(T^*,\infty)$ such that $y_n\to y^*$ locally uniformly in $(T^*,\infty)$ and $\{y^\prime_n\}$ also converges locally uniformly in $(T^*,\infty)$. Furthermore, $y_n'\to (y^*)'$ pointwise in $(T^*,\infty)$ (see {\sc rudin} \cite[Theorem~7.17]{ru}). 

For an integer $m\geq 0$ and any $n$ such that $\gamma_n>m+s_0$, we define $t_m(\gamma_n)$ to be the point at which $y_n(t_m(\gamma_n))=m+s_0$. We claim that $S_m\eqdef \displaystyle\limsup_{n \to \infty}t_m(\gamma_n)<\infty$. To see this,  define $z_n(t)=y_n(t)-m-s_0$ for $\gamma_n>m+s_0$. Then, $z_n$ solves the equation
\begin{equation*}
-(\vert z_n'\vert^{N-2}z_n')'=e^{-t}f_0(z_n+m+s_0)\eqdef e^{-t}\bar{f}_0(z_n).
\end{equation*}
Let $T_m(\gamma_n)$ be the first zero of $z_n(t)$ as $t$ decreases from infinity. It can be checked that  $\bar{f}_0$ 
also satisfies assumptions (A1)-(A5) in \cite{GiPrSr-asy}. Therefore, again by \cite[Proposition~4.2]{GiPrSr-asy}, we get $S_m<\infty$.

Fix $m$. If necessary, by restricting to a further subsequence of $\{y_n\}$ (depending on $m$) so that $\displaystyle\lim_{n\to\infty}t_m(\gamma_n)=S_m$, we obtain
\begin{equation*}
y^*(S_m)=\displaystyle\lim_{n\to\infty}y_n(S_m)=\displaystyle\lim_{n\to\infty}y_n(t_m(\gamma_n))=m+s_0.
\end{equation*}
Since $y^*$ is nondecreasing, we obtain that $y^*(t)\to\infty$ as $t \to \infty$. Integrating the O.D.E. satsfied by $y_n$ we find
\begin{equation*}
(y_n')^{N-1}(s)-(y_n')^{N-1}(t)=\int_s^t f_0(y_n(\rho))e^{-\rho}\,\mathrm{d}\rho,\quad T^*<s<t<\infty.
\end{equation*}
Using the convergence of $y_n$, $y_n'$, we can pass to the limit as $n\to\infty$ on either side of the above equation to obtain that $y^*$ also satisfies the same integral equation. That is, $y^*$ solves the equation
\begin{equation*}
-(((y^*)^\prime)^{N-1})^\prime=e^{-t}f_0(y^*)\quad\mbox{ in }(T^*,\infty).
\end{equation*}
Now, from \eqref{eq3.2},
\begin{equation}\label{eq3.3}
\displaystyle\sup_{n}\int_{T(\gamma_n)}^\infty f_0(y_n)e^{-t}\,\mathrm{d}t<\infty.
\end{equation}
We now come to the value $y^*$ at $T^*$. Since $y^*$ is nondecreasing, $y^*$ has a right limit at $T^*$. Integrating the ODE satisfied by $y_n$ between $t\in [T(\gamma_n),t_0(\gamma_n)]$ and $t_0(\gamma_n)$ and using \eqref{eq3.1} and \eqref{eq3.3}, we deduce that $\{y_n'\}$ is uniformly bounded in $[T(\gamma_n), T^*+1]$. Consequently, the extended sequence $\{y_n\}$ is uniformly bounded in the Lipschitz norm on $[T^*-1, T^*]$. Then, by Ascoli-Arzela theorem,
 we have
\begin{equation*}
y^*(T^*)=\displaystyle\lim_{n\to\infty}y_n(T^*)=0.
\end{equation*}
From \eqref{eq3.3} and Fatou's lemma, we obtain that $\int_{T^*}^\infty f_0(y^*)e^{-t}\,\mathrm{d}t<\infty$. Thus, to summarize, $y^*$ solves the problem \eqref{ind} with $T=T^*$ with the additional property: $y(T^*)=0$. 

Let $R=Ne^{-T^*/N}$. We now define $z^*(x)=y^*\big(N\log(N/\vert x\vert)\big)$ for $x\in B_{R}\setminus\{0\}$.  It follows that $z^*$ solves the following problem:
$$\hspace{25mm}\left\{
\begin{array}{cllll}\left.
\begin{array}{rllll}
 -\Delta_N z^* & = & f_0(z^*)\\
z^* &>& 0 \end{array}\right\} \;\; \text{in} \; B_{R}\setminus\{0\}, \\
z^*=0 \;\; \text{on} \;\; \partial B_{R}, \;\;
\lim_{|x| \to 0} z
^*(x)  = \infty , \;\;\;
f_0(z^*) \in L^1(B_{R}).
\end{array}\right.
$$
(ii) We note that  $f$ satisfies assumptions $(f_1)-(f_3)$. From Theorem \ref{1}, $z^*$ satisfies the
 equation $-\Delta_N z^*= f(z^*)+\alpha \delta_0$ in the sense of distributions
 in  $B_{\frac{R}{2}}$ for some $\alpha \geq 0.$ Since $f$ is super-exponential,  by Lemma~\ref{13} we must have $\alpha=0$. Thus,  $z^*$ is the required singular solution  for $(P_{0, \frac{R}{2}})$.

If $z^* \in W^{1,N}_{loc}(B_R)$, by Trudinger-Moser imbedding \cite{Mo}, we obtain that $f(z^* )\in L^p_{loc}(B_R)$ for all $p \geq 1$ and hence $z^*$ is locally bounded in $B_R$, a contradiction. 
\hfill \qed
%
\section{Asymptotic behaviour in the super-exponential case }\label{section4}
%
%
From the assumption $(f_3)$, we can  fix $\kappa \geq 0$ so that the map 
\begin{equation}\label{spi}
t \mapsto F(t):= f(t) + \kappa t^{N-1} \;\text{ is strictly increasing for } \; t \geq 0.
\end{equation}

For instance, if $f(t)= e^{{t}^\mu}$, $\mu>0$, then we may choose $\kappa =0$ and check that $F^{-1}(t) = (\log t)^{\frac{1}{\mu}}$. We have the following asymptotic estimate for $u$.
\begin{lemma}\label{19} Let $f$ be a super-exponential nonlinearity and $\kappa$ be chosen as in \eqref{spi}. Assume that
\begin{equation}\label{saj}
(f(t))^\lambda \leq c f(\lambda t) \;\text{ for some } \; c>0 \; \text { and all }\; \lambda >1, \; t \geq 0.
\end{equation}
Then, any unbounded solution $u$ of $(P^*_R)$ satisfies the following properties: \\

\noindent (a) Given $\epsilon>0$, there exists $R_\epsilon \in (0, R)$ such that  $u(x) \leq F^{-1}( \epsilon |x|^{-N})$ for all $x\in B_{R_\epsilon}$.

\noindent (b)  For any  $\epsilon>0,$
$$   \limsup_{x \to 0} \frac{u(x)}{F^{-1}( \epsilon |x|^{-N})}  =1. $$ 
\end{lemma}

%

\noindent \proof We write the problem $(P^*_R)$ in radial co-ordinates as follows (see Proposition \ref{p2}):
 \begin{eqnarray*}
\left\{\begin{array}
{ll}
& -r^{1-N}\Big(\Big( r^{N-1} |u^\prime|^{N-2} u^{\prime}  \Big)\Big)^\prime
 =  f(u) \;\mbox{ in }\, {\mathcal D}^\prime((0,R)),\\
&u \;\text{ nonnegative in }\; (0,R) .\\
\end{array}\right.
 \end{eqnarray*}
 We use again the Emden-Fowler transformations as in \eqref{17}-\eqref{18} and denote by $y$ the transformed solution corresponding to $u$.  
 
Since $f$ is super-exponential, as in the proof of Lemma  \ref{13}, one can show that
\begin{equation} \label{19.3}
\lim_{t \to \infty} y^{\prime}(t)=0.
\end{equation}
Since $y$ is a concave function, it follows that $y^{\prime}(t)>0$ for all $t$.
Integrating the ODE satisfied by $y$ (see \eqref{18}) between the limits $1<< t <t_1$, we obtain
\begin{eqnarray*}
\big( y^\prime \big)^{N-1}(t) - \big( y^\prime \big)^{N-1}(t_1) &=& \int_{t}^{t_1} \big(f(y(s))+\kappa s^{N-1} \big)e^{-s} ds - \int_{t}^{t_1} \kappa s^{N-1}e^{-s} ds \notag \\
&\geq & \big(f(y(t))  + \kappa t^{N-1} \big)\big(e^{-t}-e^{-t_1}\big) - \int_{t}^{t_1} \kappa s^{N-1}e^{-s} ds. \notag \\
\end{eqnarray*}
Letting $t_1 \to \infty$ in the above inequality and using \eqref{19.3}, we deduce
\begin{equation} \label{19.4}
\big( y^\prime \big)^{N-1}(t) \geq \big(f(y(t))  + \kappa t^{N-1} \big)e^{-t}- \int_{t}^{\infty} \kappa s^{N-1}e^{-s} ds \geq - \int_{t}^{\infty} \kappa s^{N-1}e^{-s} ds.
\end{equation}
In particular, using \eqref{19.3} again, we have
$$
\lim_{t \to \infty}  f(y(t))  e^{-t}  =0.
$$
Hence, from assumptions (ii) and (iii) on $f$, we obtain that
$$
\lim_{t \to \infty} F(y(t))e^{-t} = \lim_{t \to \infty}  \big( f(y(t)) + \kappa y^{N-1}(t) \big) e^{-t}  =0.
$$
Now, (a) follows easily. Using part (a), we have that 
\begin{equation}\label{sup1}
\limsup_{x \to 0} \frac{u(x)}{F^{-1}( \epsilon |x|^{-N})}  \leq 1. 
\end{equation}  
Suppose for some $\epsilon>0$ the above inequality is strict. Then we can find $0< \eta <1$ and $R_0>0$ such that 
$$
u(x) \leq \eta F^{-1}(\epsilon  |x|^{ -N}) \;\;\;  \mbox{ in } B_{R_0}\setminus\{0\}.
$$
Using the assumption \eqref{saj} we easily obtain the pointwise estimate
$$
f(u) \leq O(1) |x|^{-\eta N} \;\;\;  \mbox{ in } B_{R_0}\setminus\{0\}.
$$
Thus,  $f(u) \in L^p_{loc}(B_R)$ for some $p>1$. From classical estimates in Serrin \cite{Se1} we deduce $u \in L^{\infty}_{loc}(B_R)$, a contradiction.  This shows that for any $\epsilon>0$ we must have equality in \eqref{sup1} which completes our proof.
\hfill \qed

\begin{cor}\label{hon}
Let $f, u$ be as in Lemma \ref{19}.  Then, one of the following  alternatives holds : \\

\noindent (i) $\displaystyle \lim_{x \to 0} \frac{u(x)}{F^{-1}( |x|^{-N})} =1$, or\\

\noindent (ii) for some $0 < c_* < 1$ and any $a \in (c_*, 1)$,  the graph of $u$ crosses the  graph of the function  $a F^{-1}( |x|^{-N})$ infinitely often in $B_R$. 

\end{cor}

\proof Let 
$$
c_* :=  \liminf_{x \to 0}\frac{u(x)}{F^{-1}( |x|^{-N})} \; .
$$
If $c_* =1$, the first alternative holds. Otherwise, we may take any $a \in (c_*,1)$ and verify that the second alternative holds for all such $a$.
\hfill\qed

We point out that a similar alternative holds in the super-critical case in higher dimensions (see Theorem 5.13 in \cite{Ve1}).  See   Corollary  \ref{meh} for a more precise version.

\section{An integral-type lower bound for singular solutions}

\begin{lemma}\label{vic}
Let $f$  satisfy the assumptions in $\eqref{spi}$  and \eqref{saj}. Further assume that there exists a $\mu>1$ such that 
\begin{equation} \label{esta}
f(t) \geq C e^{t^\mu} \;\mbox{ for all large } \; t>0 \;\mbox{ and some }\; C>0.
\end{equation}
If $u$ is an unbounded solution of $(P^*_{R})$, then the following holds:
\begin{equation}\label{bob}
 \int_0^{R} \frac{1}{|r (\log(N/r))|^{\frac{1}{\mu}}} \Bigg(  \int_{B_r(0)} f(u) dx \Bigg)^{\frac{1}{N-1}} dr = \infty .\end{equation}
\end{lemma}

\begin{remark}
Note that we always have $\int_{B_R}  f(u) dx < \infty.$
\end{remark}

\noindent \proof  Choose $\kappa>0$ as in \eqref{spi} and for any $t\geq 0$ define $G_\kappa (t):= C e^{t^\mu}+\kappa t$. Since 
$$
F(t) \geq G_\kappa(t) \geq G_0(t),  \; \mbox{ for all } t \geq 0,
$$
we have 
\begin{equation}\label{sur}
(\log t )^{\frac{1}{\mu}} = G_0^{-1}(t) \geq F^{-1}(t)  \;\text{ for all large } \; t >0.
\end{equation}
Supposing that the conclusion of lemma does not hold and write the solution $u$  in the form:
\begin{equation}\label{gin}
u(x):= v(x) (N \log(N/|x|))^{\frac{1}{\mu}} \quad \text{ for } \; x \in B_R\setminus\{0\}.
\end{equation}
From \eqref{sur} and Lemma \ref{19} (a), we know that
\begin{equation}\label{gia}
0< v(x) < 1  \;\text{ for all small} \; |x|.
\end{equation}
Going to the radial variable $r$ and the Emden-Fowler variable $t$ as given in \eqref{17}, we have that $y(t):=u(|x|)$ and $z(t):= v(|x|)$ satisfy:
\begin{equation}\label{gir}
y(t):= z(t) t^{\frac{1}{\mu}} \;\text{ for } \; t \geq T:= N \log(N/R).
\end{equation}
Since $f$ is super-exponential, we have that (see the proof of Lemma \ref{13})
$$
\lim_{t \to \infty} y^{\prime}(t)=0.
$$
Recall that $y$ solves the differential equation in \eqref{18}. We see then that $y$ also solves the equation 
\begin{equation}\label{git}
y^{\prime} (t)= \Big(\int_t^{\infty} f(y(s)) e^{-s} ds \Big)^{\frac{1}{N-1}}, \;\text{ for } \; t > T.
\end{equation}
Thus, we obtain the following equation for $z$:
\begin{equation}\label{gil}
z^{\prime}(t) + \frac{1}{\mu t} z(t)  =   t^{-\frac{1}{\mu}} \Big(\int_t^{\infty} f(z(s) s^{\frac{1}{\mu}}) e^{-s} ds \Big)^{\frac{1}{N-1}}, \; \text{ for } \; t > T.
\end{equation}
Note that $0< z(t) <1$ for all large $t>0$. 
Making the final transformation
\begin{equation}\label{gim}
\xi := \log t, \; \rho(\xi)=z(t), \;\text{ for } \; t \geq T_0:=  \max\{1, T\},
\end{equation}
the  equation in \eqref{gil} simplifies to:
\begin{equation}\label{gips}
\rho^{\prime}(\xi) + \frac{1}{\mu} \rho (\xi)  = H(\xi), \; \;\;\; \xi > \log T_0
\end{equation}
where
\begin{equation}\label{gla}
H(\xi) :=   e^{\big(\frac{\mu -1}{\mu} \big) \xi} \Big(\int_{\xi}^{\infty} f(\rho(\zeta) e^{\frac{\zeta}{\mu}}) e^{-e^\zeta} e^{\zeta}d\zeta \Big)^{\frac{1}{N-1}}.
\end{equation}
Note that $0< \rho(\xi) <1$ for all large $\xi>0$. It can also be deduced, following the proof of Lemma \ref{19}(b),  that 
$$
 \limsup_{\xi \to \infty} \rho(\xi) =1.
$$
The assumption that the integral in  \eqref{bob} is finite translates  to
\begin{equation}\label{as}
H_0:= \int_{\log T_0}^{\infty} H(\xi) d\xi < \infty. 
\end{equation}

Integrating the equation \eqref{gips} between the limits $\xi > \log T_0$ and $\xi_n > \xi$, we obtain :
\begin{equation}\label{joci}
\rho(\xi_n)- \rho(\xi)+  \frac{1}{\mu} \int_{\xi}^{\xi_n} \rho(\zeta)d\zeta = \int_{\xi}^{\xi_n} H(\zeta)d\zeta   \;\;\;\;  \text{ for all } \; \xi > \log T_0.
\end{equation}
 Letting $n \to \infty$, using the fact that $0 < \rho <1$ and \eqref{as},
 \begin{equation}\label{joc}
\frac{1}{\mu} \int_{\xi}^{\infty} \rho(\zeta)d\zeta \leq \int_{\xi}^{\infty} H(\zeta)d\zeta + 1 \leq H_0 +1  \; \text{ for all } \; \xi > \log T_0.
\end{equation}
Combining with \eqref{gips},  we get that  $\rho \in W^{1,1}([\log T_0, \infty))$.
 Necessarily,
$\rho(s) \to 0$ as $s \to \infty$, a contradiction to the fact that $\limsup_{s \to \infty} \rho =1$. This proves the lemma. \hfill \qed 

\begin{remark}\label{tel}
The example in lemma \ref{mu-lemma} is relevant in the context of the above result.
\end{remark}

\begin{remark} When $N=2$, by Fubini theorem, the condition \eqref{bob} reduces to:
\begin{equation}
\int_{B_R} | \log(N/|x|)|^{1-\frac{1}{\mu}}  f(u) dx  = \infty. \notag
\end{equation}
\end{remark}

Define
$$
w(r) := \int_r^{R} s^{1-N} |\log(N/s)|^{\frac{1-N}{\mu}} ds.
$$

\begin{cor}\label{coo}
Let $f$  satisfy the assumptions in Lemma \ref{vic}.  If $u$ is an unbounded solution of $(P^*_{R})$, then the following holds:
\begin{equation}\label{bobi}
\int_{B_R} f(u) w(|x|) dx = \infty. 
\end{equation}
\end{cor}
\proof
Let $R_0 := \min\{N, R\}$.  By \eqref{bob}, Jensen's inequality and Fubini theorem,
\begin{eqnarray*}
\infty &=& \Bigg (\int_0^{R_0} \frac{1}{r (\log(N/r))^{\frac{1}{\mu}}} \Big(  \int_{B_r(0)} f(u) dx \Big)^{\frac{1}{N-1}} dr \Bigg)^{N-1} \notag \\ 
&\leq& R_0^{N-2}  \int_0^{R_0} \frac{1}{r^{N-1} (\log(N/r))^{\frac{N-1}{\mu}}}  \Big(\int_{B_r(0)} f(u) dx \Big) dr \notag \\
&\leq& R_0^{N-2}  \int_{B_{R_0}} f(u) w(|x|) dx.
\end{eqnarray*}

\begin{cor}\label{bb}
Let $f(t):= e^{t^{\mu}}$, $\mu>1,$  and $u$ be a solution of $(P^*_{R})$. Suppose that for some $\theta \geq (1-\frac{1}{\mu})(N-1)+1$, 
\begin{equation}\label{wh}
u(x) \leq \Big ( N \log (N/|x|) - \theta  \log(N\log (N/|x|)) \Big)^{\frac{1}{\mu}}, \; \text{ for all small } \; |x|.
\end{equation}
Then, $u$ is bounded.
\end{cor}
\proof  It is easy to check that \eqref{wh} implies 
(refer to \eqref{gla}),  
\begin{equation}\label{geo}
\int_0^{\infty} H(\zeta)d\zeta < \infty.
\end{equation}
From proof of Lemma \ref{vic} we obtain that $\lim_{t \to \infty} z(t)=0$. Following the arguments at the end of proof of Lemma \ref{19}, the conclusion follows. \hfill \qed  \\

We can now refine the bound in Lemma \ref{19}(a) and Corollary \ref{hon} (ii)  for the nonlinearity in the above result.

\begin{cor}\label{meh}
Let $f(t):= e^{t^{\mu}}$, $\mu>1,$ and $\theta \geq (1-\frac{1}{\mu})(N-1)+1$. Then, one of the following  alternatives holds for an unbounded solution $u$ of $(P^*_{R})$: \\

\noindent (i) $\displaystyle \Big ( N \log (N/|x|) - \theta  \log(N\log (N/|x|)) \Big)^{\frac{1}{\mu}} \leq u(x) \leq \Big ( N \log (N/|x|)  \Big)^{\frac{1}{\mu}}, \; \text{ for all small } \; |x|$ \quad or\\

\noindent (ii) $u(x) \leq \Big ( N \log (N/|x|)  \Big)^{\frac{1}{\mu}} \; \text{ for all small } \; |x|$ and  the graph of $u$ crosses the  graph of the function  $\Big ( N \log (N/|x|) - \theta  \log(N\log (N/|x|)) \Big)^{\frac{1}{\mu}}$ infinitely often  in $B_r$ for all small $r>0$. 
\end{cor}

\end{document}